\documentclass[conference]{IEEEtran}
\usepackage[T1]{fontenc}
\usepackage{cite}
\usepackage{amsmath,amssymb,amsfonts, amsthm}
\usepackage{graphicx}
\usepackage{booktabs}
\usepackage{array}
\usepackage{multirow}
\usepackage{hyperref}
\usepackage{siunitx}
\usepackage{xcolor}
\usepackage{float}
\usepackage{caption}
\usepackage{subcaption}
\usepackage{afterpage}
\usepackage{textcomp}
\usepackage{caption}
\usepackage{mathtools}
\usepackage{graphicx}
\usepackage{cite}
\usepackage{algorithm}
\usepackage{algorithmic}
\usepackage{multirow}

\newcommand{\N}{\mathbb{N}}

\theoremstyle{definition}
\newtheorem{definition}{Definition}
\newtheorem{example}{Example}
\newtheorem{remark}{Remark}

\def\BibTeX{{\rm B\kern-.05em{\sc i\kern-.025em b}\kern-.08em
    T\kern-.1667em\lower.7ex\hbox{E}\kern-.125emX}}

\usepackage{fancyhdr}
\begin{document}

\title{Frequency Sensitive Duplicate Detection Using Multi-Metric Spaces}

\author{
\textbf{Debjyoti Chatterjee} \\
Jalpaiguri Government Engineering College, Jalpaiguri, West Bengal, India  \\
debjyoti.chatterjee@jgec.ac.in
\\[0.8em]

\textbf{Shashi bajaj Mukherjee} \\
Dr. B. C. Roy Engineering College, Durgapur, West Bengal, India \\
shashiagarwala@gmail.com
\\[0.8em]

}

\maketitle
\thispagestyle{fancy}

\begin{abstract}
Classical metric spaces often fail to model data-intensive systems where repetition and frequency of values are meaningful. In applications such as transactional databases, sensor logs, and record linkage, conventional distance measures ignore multiplicity information, leading to information loss and incorrect similarity judgments. This paper introduces multi-metric spaces defined on multisets and valued in the multi-real number system, providing a principled way to incorporate frequency into distance computations. We demonstrate the usefulness of multi-metrics through a frequency sensitive duplicate detection example, showing improved accuracy over classical metric based approaches.
\end{abstract}


\section{Introduction}
\label{sec:intro}

Modern engineering systems and data driven applications routinely generate records that contain repeated observations, aggregated measurements, or frequency based summaries rather than simple collections of distinct values. Examples include transactional databases where items may be purchased multiple times, sensor networks that record repeated readings over fixed intervals, and log or clickstream data that encode user behavior through counts and intensities. In such settings, repetition is not accidental noise but a meaningful indicator of importance, strength, or frequency of occurrence, and any analytical framework that ignores this information risks significant semantic distortion.

Despite this reality, most classical similarity and distance models are built on metric spaces defined over sets or fixed length vectors, where repeated elements are either collapsed into single occurrences or normalized away. As a consequence, records that differ substantially in their frequency distributions may be judged as identical or nearly identical, leading to misleading similarity assessments. This limitation is particularly problematic in tasks such as duplicate detection, clustering, and data integration, where accurate discrimination between records depends not only on which attributes are present but also on how often they occur.

To overcome these shortcomings, we adopt the framework of multi-metric spaces \cite{ChatterjeeMukherjee2025b} defined on multisets, in which distances are valued in the multi-real number system \cite{ChatterjeeMukherjee2025b,ChatterjeeSamanta2020, ChatterjeeSamanta2018}. This approach extends classical metric spaces by explicitly incorporating multiplicity into the distance computation, thereby preserving both attribute identity and frequency information. By doing so, multi-metric spaces provide a mathematically rigorous and application oriented foundation for analyzing frequency sensitive data, making them well suited for modern problems in data cleaning, pattern recognition, clustering, and decision making within data intensive systems.

\section{Preliminaries}\label{}

\subsection{Multisets: Definitions and Basic Operations}\label{subsec2}
The concept of multisets (abbreviated as msets) was first introduced by Yager~\cite{Yager1986}. Fundamental definitions, along with the notions of relations and functions in the multiset framework, were subsequently developed by Girish and John~\cite{GirishJacob2011,GirishJacob2012a,GirishJacob2012b}.

An mset \cite{GirishJacob2011,GirishJacob2012a,GirishJacob2012b} M drawn from the set X is represented by a function $Count_M$ or $C_{M}$ defined as $C_{M}:X\rightarrow \mathbb{N}$, Where $N$ is the set of all non-negative integers. $C_{M}(x)$ is the number of occurrences of the element $x\in X$ in the mset $M$.
Let $M$ a mset drawn from the set $X=\{x_{1}, x_{2}, ..., x_{n}\}$ with $x_{i}$ appearing $k_{i}$ times in M, then it is denoted by $x_{i}\in ^{k_{i}} M$. The mset $M$ drawn from the set $X$ is then denoted by $\{k_{1}/x_{1},k_{2}/x_{2},...,k_{n}/x_{n}\}$. However, those elements that are not included in the mset $M$ have zero count. Clearly, a crisp set is a special case of an mset. Let $X=\{a,b,c,d,e\}$ be any set. Then $M=\{3/a,2/b,1/e\}$ is an mset drawn from $X$.\\
Let $P$ and $Q$ be two msets drawn from a set $X$, then the following are defined \cite{GirishJacob2011,GirishJacob2012a,GirishJacob2012b}:\\
\
\noindent (i) $P=Q$ if $C_{P}(x)=C_{Q}(x)$ $\forall x\in X$.\\
\noindent (ii) $P\subseteq Q$ if $C_{P}(x)\leq C_{Q}(x)$ $\forall x\in X$, then we call $P$ to be a \textbf{submultiset} (submset, in short) of $Q$.\\
\noindent (iii) $M=P\cup Q$ if $C_{M}(x)=max\{C_{P}(x),C_{Q}(x)\}$ $\forall x\in X$.\\
\noindent (iv) $M=P\cap Q$ if $C_{M}(x)=min\{C_{P}(x),C_{Q}(x)\}$ $\forall x\in X$.\\
\noindent (v) $M=P\oplus Q$ if $C_{M}(x)=C_{P}(x)+C_{Q}(x)$ $\forall x\in X$.\\
\noindent (vi) $M=Q\ominus P$ if $C_{M}(x)=max\{C_{Q}(x)-C_{P}(x),0\}$ $\forall x\in X$.\\
Here $\cup$, $\cap$, $\oplus$ and $\ominus$ represent mset union, mset intersection, mset addition and mset subtraction, respectively.\\

Let M be an mset drawn from a set $X$, then the \textbf{support set} of $M$ denoted by $M^{*}$ is a subset of $X$ and $M^{*}=\{x\in X:C_{M}(x)>0\}$. i.e., $M^{*}$ is an ordinary set and it is also called root set. An mset $M$ is is said to be an empty multiset if $C_{M}(x)=0$ for all $x\in X$. The \textbf{cardinality} of an mset $M$ drawn from a set $X$ is denoted by $card(M)$ or $|M|$ and is given by $|M|=\sum_{x\in X}C_{M}(x)$.
The \textbf{bounded mset space} $[X]^{w}$ is the set of all msets whose elements are in X such that no element in the mset occurs more than $w$ times. The \textbf{unbounded mset space} $[X]^{\infty}$ is the set of all msets drawn from $X$ such that no limit on the number of occurrences of an object in an mset.\\
Let $\{M_{i}:i\in \Omega\}$ be a collection of msets drawn from $[X]^{m}$, then the following operations are defined:\\
\noindent (i) $P=\bigcup_{i\in \Omega} M_{i}$ if $C_{P}(x)=max_{i\in \Omega}C_{M_{i}}(x), x\in X$.\\
\noindent (ii) $P=\bigcap_{i\in \Omega} M_{i}$ if $C_{P}(x)=min_{i\in \Omega}C_{M_{i}}(x), x\in X$.\\
Let $X$ be a non-empty set and $[X]^{w}$ be the boundrd mset space defined over $X$, then the complement $M^{c}$ of $M$ in $[X]^{w}$ is an element of $[X]^{w}$ such that \cite{GirishJacob2011,GirishJacob2012a,GirishJacob2012b} $C_{M}^{c}(x)=w-C_{M}(x)$ $\forall x\in X$.\\

\section{Multi Points}

\begin{definition}
Let $M$ be a multiset over a universal set $X$.
Let $x\in M^{*}$ and $k\in \mathbb{N}$ such that $0<k\leq C_M(x)$, then the \emph{multi point} \cite{DasRoy2021} $P_x^k$ of $M$ is the multiset over $X$ is defined as follows: \\
for $y\in X$, $P_x^k(y)=\begin{cases} k, & y=x, \\ 0, & y\neq x. \end{cases}$\\

The element $x$ is called the \emph{base} and $k$ is called the \emph{multiplicity} of the multi point $P_x^k$. Collection of all multi points of an mset $M$ is denoted by $M_{pt}$. i.e., $M_{pt} = \{P_x^k : x\in M^{*}, 0< k \le C_M(x), k\in N\}$
\end{definition}

.

\section{Multi Real Numbers}

\begin{definition}
For all $a\in R$ and for all $k\in \mathbb{N}_0$, where $\mathbb{N}_0=\{0,1,2,...\}$, let us define $R_a^k$ to be a multi real number \cite{ChatterjeeSamanta2020,ChatterjeeSamanta2018,ChatterjeeMukherjee2025b}. The set of all multi real numbers is denoted by $m(\mathbb{R})$. For a multi real number $R_a^k$, if $a>0$, then it is called a non-negative multi real number. The set of all non-negative multi real number is denoted by $m(\mathbb{R^+})$.
\end{definition}

\begin{remark}
Each multi real number \(R_a^k \in m(\mathbb{R})\) may be interpreted as a pair
\((a,k)\), where \(a\) represents the numerical value and \(k\) encodes its
multiplicity or frequency. Thus, \(m(\mathbb{R})\) extends the classical real
number system by incorporating occurrence information, while retaining the
real line as its value component.
\end{remark}

\begin{definition}
Let $R_a^i$ and $R_b^j$ be two multi real numbers. We define \cite{DasRoy2021}
\[
R_a^i = R_b^j
\]
if $a = b$ and $i = j$.

Let $R_a^i$ and $R_b^j$ be two multi real numbers. We define \cite{DasRoy2021}
\[
R_a^i < R_b^j
\]
if either $a < b$, or ($a = b$ and $i < j)$.

Let $R_a^i$ and $R_b^j$ be two multi real numbers. We define \cite{DasRoy2021}
\[
R_a^i \le R_b^j
\]
if either $R_a^i = R_b^j$, or $R_a^i < R_b^j$.\\
$\le$ is immediately a total order.
\end{definition}

\begin{definition}
Addition and multiplication of two multi real numbers are defined by\\
$$R_a^k \oplus R_b^\ell = R_{a+b}^{k+\ell}$$
$$R_a^k \otimes R_b^\ell = R_{ab}^{k\ell}$$

\end{definition}

\subsection{Theorem}
The algebraic structure $(m(\mathbb{R}), \oplus, \otimes)$ is a commutative semiring, and the order relation $\leq$ defined in Definition~3 is a total order compatible with both $\oplus$ and $\otimes$.

\begin{proof}
Let $R_a^k, R_b^\ell, R_c^m \in m(\mathbb{R})$.

\medskip
\noindent
\textbf{(i) Closure.}
By Definition~4,
\[
R_a^k \oplus R_b^\ell = R_{a+b}^{k+\ell}, \qquad
R_a^k \otimes R_b^\ell = R_{ab}^{k\ell}.
\]
Since $a+b, ab \in \mathbb{R}$ and $k+\ell, k\ell \in \mathbb{N}_0$, both results belong to $m(\mathbb{R})$. Hence, $m(\mathbb{R})$ is closed under $\oplus$ and $\otimes$.

\medskip
\noindent
\textbf{(ii) Associativity.}
\[
\begin{aligned}
(R_a^k \oplus R_b^\ell) \oplus R_c^m
&= R_{a+b}^{k+\ell} \oplus R_c^m
= R_{a+b+c}^{k+\ell+m} \\
&= R_a^k \oplus R_{b+c}^{\ell+m}
= R_a^k \oplus (R_b^\ell \oplus R_c^m),
\end{aligned}
\]
and
\[
\begin{aligned}
(R_a^k \otimes R_b^\ell) \otimes R_c^m
&= R_{ab}^{k\ell} \otimes R_c^m
= R_{(ab)c}^{(k\ell)m} \\
&= R_{a(bc)}^{k(\ell m)}
= R_a^k \otimes (R_b^\ell \otimes R_c^m).
\end{aligned}
\]

\medskip
\noindent
\textbf{(iii) Commutativity.}
\[
R_a^k \oplus R_b^\ell = R_{a+b}^{k+\ell} = R_{b+a}^{\ell+k} = R_b^\ell \oplus R_a^k,
\]
\[
R_a^k \otimes R_b^\ell = R_{ab}^{k\ell} = R_{ba}^{\ell k} = R_b^\ell \otimes R_a^k.
\]

\medskip
\noindent
\textbf{(iv) Identity elements.}
Define $0_m = R_0^0$ and $1_m = R_1^1$. Then,
\[
R_a^k \oplus R_0^0 = R_a^k, \qquad
R_a^k \otimes R_1^1 = R_a^k.
\]
Thus, $0_m$ and $1_m$ act as the additive and multiplicative identities, respectively.

\medskip
\noindent
\textbf{(v) Distributivity.}
\[
\begin{aligned}
R_a^k \otimes (R_b^\ell \oplus R_c^m)
&= R_a^k \otimes R_{b+c}^{\ell+m}
= R_{a(b+c)}^{k(\ell+m)} \\
&= R_{ab+ac}^{k\ell+km}
= R_{ab}^{k\ell} \oplus R_{ac}^{km} \\
&= (R_a^k \otimes R_b^\ell) \oplus (R_a^k \otimes R_c^m).
\end{aligned}
\]

Hence, $(m(\mathbb{R}), \oplus, \otimes)$ is a commutative semiring.

\medskip
\noindent
\textbf{(vi) Compatibility of the order $\leq$.}

Recall that
\[
R_a^k \le R_b^\ell \iff (a < b) \ \text{or} \ (a=b \ \text{and}\ k \le \ell).
\]

If $R_a^k \le R_b^\ell$, then either $a<b$, which implies $a+c < b+c$, or $a=b$ and $k \le \ell$, which implies $k+m \le \ell+m$. Hence,
\[
R_a^k \oplus R_c^m \le R_b^\ell \oplus R_c^m.
\]

Similarly, if $R_c^m \in m(\mathbb{R}^+)$, then $a<b$ implies $ac<bc$, and $k \le \ell$ implies $km \le \ell m$. Therefore,
\[
R_a^k \otimes R_c^m \le R_b^\ell \otimes R_c^m.
\]

Thus, the order $\leq$ is compatible with both $\oplus$ and $\otimes$.

\end{proof}

\section{Multi-Metric Spaces}

\begin{definition}
Let $X$ be a nonempty set. Let $M$ be a multiset drawn from $X$ having multiplicity of any element atmost equal to $w\in \mathbb{N}$. A \emph{multi-metric} $d$ on $M$ is a mapping
\[
d : M_{pt} \times M_{pt} \longrightarrow m(\mathbb{R^+}).
\]

such that for all $P^i_x,P^j_y,P^k_z\in X$:
\begin{enumerate}
\item $d(P^i_x,P^j_y)=R_0^0$ if and only if $P^i_x=P^j_y$;
\item $d(P^i_x,P^j_y)=d(P^j_y,P^i_x)$;
\item $d(P^i_x,P^k_z) \le d(P^i_x,P^j_y) \oplus d(P^j_y,P^k_z)$.
\end{enumerate}
\end{definition}
Then, \((M,d)\) is called a multi-metric space.
\begin{example}
Let
\[
X=\{a,b,c\}
\]
and let \(M\) be a multiset drawn from \(X\) such that the multiplicity of each element is at most \(5\in \mathbb{N}\).

Define a mapping
\[
d : M_{pt}\times M_{pt}\longrightarrow m(\mathbb{R}^+)
\]
by
\[
d(P^i_x,P^j_y)=
\begin{cases}
R_0^{0}, & \text{if } P^i_x=P^j_y,\\[6pt]
R_1^{\,|i-j|+1}, & \text{if } R^i_x\neq R^j_y.
\end{cases}
\]
Hence, \(d\) is a multi-metric on \(M\), and \((M,d)\) is a multi-metric space.
\end{example}

\begin{example}
Let \((X,\rho)\) be a classical metric space, where
\[
X=\mathbb{R}, \qquad \rho(x,y)=|x-y|, x,y\in \mathbb{R}.
\]
Let \(M\) be a multiset drawn from \(X\) such that each element has finite multiplicity at most 7.

Define a mapping
\[
d : M_{pt}\times M_{pt}\longrightarrow m(\mathbb{R}^+)
\]
by
\[
d(P^i_x,P^j_y)=R^{\,|i-j|}_{\rho(x,y)}.
\]

Hence, \(d\) is a multi-metric on \(M\), and \((M,d)\) is a multi-metric space.
\end{example}

\begin{remark}
The multi real number \(0_m = R_0^0\) represents the \emph{multi-zero}. It denotes
the absence of both numerical magnitude and multiplicity. In particular,
\(R_0^0\) plays the role of the neutral element for the addition operation
\(\oplus\) in \(m(\mathbb{R})\), since for any \(R_a^k \in m(\mathbb{R})\),
\[
R_a^k \oplus R_0^0 = R_a^k.
\]
Moreover, in the context of multi-metric spaces, \(0_m = R_0^0\) naturally
corresponds to zero distance, indicating that two multi points are identical.
\end{remark}

\begin{remark}
This example shows how any classical metric can be lifted to a multi-metric by incorporating multiplicity differences into the distance value, thereby capturing both numerical deviation and frequency imbalance.
\end{remark}
`

\begin{definition}[Open Ball]
Let \((M,d)\) be a multi-metric space and let \(P^i_x \in M_{pt}\).
For a non-negative multi real number \(R_r^k \in m(\mathbb{R}^+)\),
the \emph{open ball} centered at \(P^i_x\) with multi-radius \(R_r^k\) is defined as
\[
B_d(P^i_x,R_r^k)
=
\left\{\, P^j_y \in  M_{pt}\;:\; d(P^i_x,P^j_y) < R_r^k \,\right\}.
\]
\end{definition}

\begin{definition}[Open Set]
A subset \(U \subseteq M\) is said to be \emph{open} in the multi-metric space
\((M,d)\) if for every \(P^i_x \in U_{pt}\), there exists a non-negative multi real
number \(R_r^k \in m(\mathbb{R}^+)\) such that
\[
B_d(P^i_x,R_r^k) \subseteq U_{pt}.
\]
The family of all open sets induced by the multi-metric \(d\) is denoted by \(\tau_d\).
\end{definition}

\subsection{Theorem}
Let \((M,d)\) be a multi-metric space.  
The family \(\tau_d\) of all open sets induced by the multi-metric \(d\) forms a
topology on \(M\), called the \emph{multi-metric topology}.

\begin{proof}
Let \(\tau_d\) denote the collection of all open subsets of \(M\) as defined in
Definition~7. We verify the axioms of a topology.

\medskip
\noindent
\textbf{(T1) \(\emptyset\) and \(M\) are open.}

The empty set \(\emptyset\) is open vacuously, since it contains no points.

To show that \(M\) is open, let \(P_x^i \in M_{pt}\).  
Choose any non-negative multi real number \(R_r^k \in m(\mathbb{R}^+)\).
Then
\[
B_d(P_x^i,R_r^k) \subseteq M_{pt},
\]
which implies that \(M\) is open.

\medskip
\noindent
\textbf{(T2) Arbitrary unions of open sets are open.}

Let \(\{U_\alpha\}_{\alpha \in \Lambda}\) be an arbitrary family of open sets in
\(\tau_d\), and let
\[
U = \bigcup_{\alpha \in \Lambda} U_\alpha.
\]
Take any \(P_x^i \in U_{pt}\). Then there exists some \(\alpha_0 \in \Lambda\) such
that \(P_x^i \in (U_{\alpha_0})_{pt}\).
Since \(U_{\alpha_0}\) is open, there exists a non-negative multi real number
\(R_r^k \in m(\mathbb{R}^+)\) such that
\[
B_d(P_x^i,R_r^k) \subseteq (U_{\alpha_0})_{pt} \subseteq U_{pt}.
\]
Hence, \(U\) is open.

\medskip
\noindent
\textbf{(T3) Finite intersections of open sets are open.}

Let \(U_1, U_2, \dots, U_n\) be open sets in \(\tau_d\), and let
\[
U = \bigcap_{m=1}^n U_m.
\]
Take any \(P_x^i \in U_{pt}\). Then \(P_x^i \in (U_m)_{pt}\) for each \(m=1,2,\dots,n\).
Since each \((U_m)_{pt}\) is open, there exists a non-negative multi real number
\(R_{r_m}^{k_m} \in m(\mathbb{R}^+)\) such that
\[
B_d(P_x^i,R_{r_m}^{k_m}) \subseteq (U_m)_{pt}.
\]

Let
\[
P_r^k = \min \{ R_{r_1}^{k_1}, R_{r_2}^{k_2}, \dots, R_{r_n}^{k_n} \},
\]
where the minimum is taken with respect to the total order on
\(m(\mathbb{R}^+)\).
Then
\[
B_d(P_x^i,R_r^k)
\subseteq
\bigcap_{m=1}^n B_d(P_x^i,R_{r_m}^{k_m})
\subseteq
\bigcap_{m=1}^n (U_m)_{pt}
= U_{pt}.
\]
Thus, \(U\) is open.

\medskip
Since \(\tau_d\) satisfies all three topology axioms, it follows that
\(\tau_d\) is a topology on \(M\).
\end{proof}

\section{Application: Duplicate Detection}
\label{sec:duplicate_detection}

In many real world data driven applications such as transactional databases, sensor streams, text corpora, and record linkage systems, identical attribute values may occur multiple times within a single record, with such repetitions carrying important semantic meaning related to frequency, intensity, or count. Consequently, a record is not merely a collection of distinct attribute values but often encodes rich multiplicity information that is lost under classical set based or vector based distance models, which collapse repeated elements into single instances. This limitation is particularly critical in duplicate detection tasks arising in data cleaning, database integration, and record linkage, where ignoring frequency information can lead to erroneous identification of duplicates. Multi-metric spaces address this shortcoming by providing a mathematically rigorous framework that explicitly incorporates multiplicity into distance computations, thereby enabling frequency sensitive discrimination between records that may appear identical under classical metrics but differ substantially in their underlying occurrence structures.

\subsection{Limitations of Conventional Distance Models and Record Representation as Multisets}

Traditional distance measures, such as the Hamming distance, Jaccard distance, cosine similarity, or Euclidean metric, typically operate on sets or fixed length vectors. When applied to data containing repeated attribute values, these models either collapse multiple occurrences into a single presence indicator or normalize away frequency information. As a result, crucial structural information is lost.

For instance, consider two transactional records:
\[
T_1 = \{a,a,a,b\}, \qquad T_2 = \{a,b\}.
\]
Under a set-based representation, both records reduce to the same set $\{a,b\}$, yielding zero distance, despite the clear difference in attribute multiplicities. This phenomenon leads to ambiguity and misclassification in tasks such as duplicate detection, anomaly identification, and similarity based retrieval.

To address this limitation, records can be modeled as multisets over a finite attribute universe $A$. Each record $T$ is represented by a multiplicity function
\[
C_T : A \to \mathbb{N},
\]
where $C_T(a)$ denotes the number of occurrences of attribute $a$ in $T$.

Thus, the records $T_1$ and $T_2$ above are represented as:
\[
C_{T_1}(a)=3,\; C_{T_1}(b)=1, \qquad
C_{T_2}(a)=1,\; C_{T_2}(b)=1.
\]

Let $A = \{a_1, a_2, \ldots, a_n\}$ denote a finite attribute universe.
Each database record $T$ is modeled as a multiplicity function 
\[
C_T : A \to \mathbb{N},
\]
or equivalently as a multiset
\[
T = \{ a_1^{k_1}, a_2^{k_2}, \ldots, a_n^{k_n} \},
\]
where $k_i \in \N$ represents the frequency or count of attribute $a_i$ in the record.
Such representations naturally arise in transactional data, clickstream summaries,
log files, and purchase histories.

Let $T_{pt}$ denote the collection of all multi points of $T$.
Distances between records are evaluated through a multi-metric
\[
d : T_{pt} \times T_{pt} \longrightarrow m(\mathbb{R^+}).
\]

\subsection{Multi-Metric-Based Duplication Criterion}

Let $T_i$ and $T_j$ be two candidate records represented as submultisets of a multiset $T$ over a finite attribute universe $A$.
Let the elementary distance is given by a multi-metric $d$ on $T$ of the form
\[
d\bigl(P_a^k,P_a^{\ell}\bigr) = R^{f(k,\ell)}_{0},
\]
where $f : \mathbb{N} \times \mathbb{N} \to \mathbb{N}$ is a function capturing multiplicity imbalance, such as
\[
f(k,\ell) = |k-\ell| .
\]
The multi-metric distance between $T_i$ and $T_j$ is defined as
\[
\delta(T_i,T_j) = \sum_{a \in A} d\bigl(P_a^{C_{T_i}(a)},P_a^{C_{T_j}(a)}\bigr).
\]
\begin{definition}[Duplicate Detection Criterion]
Two records $T_i$ and $T_j$ are declared duplicates if
\[
\delta(T_i, T_j) < \varepsilon,
\]
where $\varepsilon >R^0_0$ is a user-defined tolerance threshold reflecting
acceptable deviations in attribute frequencies.
\end{definition}

This criterion ensures that duplication is evaluated not only on attribute equality
but also on the similarity of their frequency distributions.

\section{Numerical Example: Multi-Metric-Based Duplication Criterion}

Let $T_i$ and $T_j$ be two candidate records represented as multisets with attribute set $A = \{a_1, a_2, a_3\}$: 

\[
T_i = \{a_1^2, a_2^1, a_3^3\}, \quad 
T_j = \{a_1^1, a_2^2, a_3^2\}
\]

\subsection*{Step 1: Compute Pairwise Distances}
Then,
\[
\begin{aligned}
\delta(T_i,T_j)
&= R^{|2-1|}_0 + R^{|1-2|}_0 + R^{|3-1|}_0 \\
&= R^1_0 + R^1_0 + R^2_0\\
&= R^4_0.
\end{aligned}
\]

\subsection*{Step 3: Duplication Decision}

Suppose the duplication threshold is $\varepsilon = R^2_0$. Then:

\[
\delta(T_i, T_j) = R^4_0 > \varepsilon
\]

Thus, $T_i$ and $T_j$ are \textbf{not duplicates} according to the multi-metric criterion.

If instead $T_j = \{a_1^2, a_2^1, a_3^3\}$, then:
\[
\begin{aligned}
\delta(T_i,T_j)
&= R^{|2-2|}_0 + R^{|1-1|}_0 + R^{|3-3|}_0 \\
&= R^0_0 + R^0_0 + R^0_0\\
&= R^0_0< \varepsilon.
\end{aligned}
\]

so $T_i$ and $T_j$ are \textbf{duplicates}.

\subsection{Analytical Example}

Consider a retail transaction database in which each record represents a customer's
purchase summary over a fixed time window. Let the attribute set be
\[
A = \{\text{Bread}, \text{Milk}, \text{Eggs}\}.
\]

Suppose two transaction records are given by
\[
T_1 = \{\text{Bread}^2, \text{Milk}^1, \text{Eggs}^3\}, \quad
T_2 = \{\text{Bread}^1, \text{Milk}^1, \text{Eggs}^3\}.
\]

\paragraph{Classical Metric Analysis.}
If records are treated as ordinary sets, both reduce to
\[
\{\text{Bread}, \text{Milk}, \text{Eggs}\},
\]
and any classical metric $d_c$ based on attribute equality yields
\[
d_c(T_1, T_2) = 0,
\]
thereby incorrectly classifying the records as duplicates.

\paragraph{Multi-Metric Analysis.}
Using a multi-metric with multiplicity-sensitive mapping $f(k,\ell) = |k - \ell|$, we obtain
\[
\delta(T_1, T_2)
= R^{|2-1|}_{0} + R^{|1-1|}_{0} + R^{|3-3|}_{0}
= R^{1}_{0}.
\]
Thus,
\[
\delta(T_1, T_2) > R^0_0,
\]
correctly indicating that the records differ due to variation in purchase frequency
of \emph{Bread}. For any threshold $\varepsilon < R^{1}_{0}$, the records are
\emph{not} classified as duplicates.

\begin{remark}
This example demonstrates that multi-metric based duplicate detection preserves
frequency level distinctions that are lost under classical metrics, preventing
false positives in data cleaning and record linkage tasks.
\end{remark}

\subsection{Methodological Advantages}

The proposed multi-metric duplicate detection framework offers several advantages:
\begin{itemize}
    \item It distinguishes records with identical attribute sets but differing
    frequency distributions.
    \item It supports tunable sensitivity through the choice of multiplicity mapping
    $f$ and threshold $\varepsilon$.
    \item It remains compatible with classical duplicate detection, which is recovered
    as a special case when all multiplicities are equal to one.
    \item It provides a topologically well-structured neighborhood system, enabling
    stable and interpretable similarity queries.
\end{itemize}

As a result, multi-metric spaces provide a robust and mathematically sound foundation
for duplicate detection in modern data intensive systems where repetition carries
semantic significance.

\subsection{Multi-Metric Duplicate Detection Algorithm}

We now present a pseudo-algorithm for duplicate detection based on multi-metric
distances. The algorithm identifies duplicate record pairs by evaluating frequency
sensitive distances between multiset represented records.

\begin{algorithm}[H]
\caption{Multi-Metric-Based Duplicate Detection}
\label{alg:duplicate_detection}
\begin{algorithmic}[1]
\REQUIRE Database $D = \{T_1, T_2, \ldots, T_n\}$ represented as multisets
\REQUIRE Multi-metric distance $\delta$
\REQUIRE Threshold $\varepsilon \in m(\mathbb{R}^+)$
\ENSURE Set $\mathcal{D}$ of duplicate record pairs

\STATE Initialize $\mathcal{D} \gets \emptyset$
\FOR{$i = 1$ to $n-1$}
    \FOR{$j = i+1$ to $n$}
        \STATE Compute $\delta(T_i, T_j)$
        \IF{$\delta(T_i, T_j) < \varepsilon$}
            \STATE $\mathcal{D} \gets \mathcal{D} \cup \{(T_i,T_j)\}$
        \ENDIF
    \ENDFOR
\ENDFOR
\RETURN $\mathcal{D}$
\end{algorithmic}
\end{algorithm}

The algorithm performs pairwise comparison of records using a frequency sensitive
multi-metric distance. Unlike classical duplicate detection methods, which collapse
multiplicities, this approach preserves frequency information and prevents false
positives caused by identical attribute sets with differing occurrence patterns.
The threshold $\varepsilon$ allows tuning of duplication sensitivity according to
application specific tolerance levels.

\subsection{Time Complexity Analysis}

Let $D = \{T_1, T_2, \ldots, T_n\}$ be a database of $n$ records, where each record
$T_i$ is represented as a multiset over a finite attribute universe
$X$ with $|X| = m$. Duplicate detection is performed by evaluating the
multi-metric distance $\delta$ between record pairs.

\begin{itemize}
\item \textbf{Number of comparisons.}
The algorithm compares all unordered pairs $(T_i,T_j)$ with $i<j$.
Hence, the total number of distance computations is
\[
\binom{n}{2} = O(n^2).
\]

\item \textbf{Cost of a single multi-metric computation.}
Evaluating $\delta(T_i,T_j)$ requires comparing the multiplicities of
attributes in the support of $T_i \cup T_j$. In the worst case, this
support coincides with $X$, yielding a computational cost of $O(m)$
per comparison.

\item \textbf{Overall time complexity.}
Combining the above, the total time complexity of the duplicate detection
algorithm is
\[
O(n^2 m).
\]
\end{itemize}

Although the quadratic dependence on the number of records is typical
for exact duplicate detection methods, the proposed approach provides
greater discriminatory power by preserving frequency information that
is ignored by classical metrics.

\subsection{Blocking and Indexing for Large Databases}

For large-scale databases, the exhaustive comparison cost $O(n^2)$
becomes computationally prohibitive. To improve scalability, standard
blocking and indexing techniques can be adapted to the multi-metric
framework while retaining multiplicity sensitivity.

\paragraph{ Blocking Based on Support Sets.}
Each record $T$ induces a support set
\[
\operatorname{supp}(T) = \{x \in X : C_T(x) > 0\}.
\]
Records with disjoint or weakly overlapping support sets cannot be
duplicates under reasonable multi-metric thresholds. Hence, the database
can be partitioned into blocks
\[
B_k = \{T_i \in D : \varphi(T_i) = k\},
\]
where $\varphi$ is a blocking key based on features such as support
cardinality, dominant attributes, or coarse multiplicity summaries.
Comparisons are then restricted to records within the same block.

\paragraph{ Multiplicity Aware Indexing.}
Unlike classical indexing schemes, multi-metric duplicate detection
requires indexing structures that preserve frequency information.
For each record $T_i$, a compact signature such as
\[
\Sigma(T_i) = \bigl(C_{T_i}(x_1), C_{T_i}(x_2), \ldots, C_{T_i}(x_m)\bigr),
\]
or a lower dimensional projection thereof, may be stored in an index.
Candidate records are retrieved using range or similarity queries that
respect admissible multiplicity deviations before exact multi-metric
evaluation.

\paragraph{ Threshold Based Pruning.}
Let $\varepsilon \in m(\mathbb{R}^+)$ be the duplication threshold.
During distance computation, partial distances over subsets
$Y \subset X$ can be used to derive lower bounds. If
\[
\delta_Y(T_i,T_j) \ge \varepsilon,
\]
then the full computation of $\delta(T_i,T_j)$ can be safely skipped.
This pruning strategy is particularly effective for sparse or
high dimensional multiset data.

\paragraph{ Optimized Complexity.}
Let $b \ll n$ denote the average block size after blocking. The number of
comparisons is reduced from $O(n^2)$ to $O(nb)$, yielding an overall
time complexity of
\[
O(n b m).
\]
In practice, when blocks are small and supports are sparse, this behavior
is close to linear in $n$.

\begin{remark}
Blocking and indexing significantly improve the scalability of
multi-metric based duplicate detection while preserving its correctness.
Unlike classical preprocessing methods, these optimizations respect
attribute multiplicities, ensuring that frequency-based distinctions are
not lost.
\end{remark}

\subsection{Streaming Duplicate Detection for Dynamic Data}

In many real world applications, data arrives continuously as a stream
(e.g., transactional logs, sensor readings, or online user records). Recomputing
all pairwise multi-metric distances upon each insertion is infeasible. A
\emph{streaming variant} of duplicate detection leverages incremental updates,
bounded memory, and approximate indexing to maintain efficient and accurate
detection.

\paragraph{Incremental Multi-Metric Maintenance.}
Let $D_t$ denote the dataset at time $t$. For a newly arriving record $T_{new}$:
\begin{enumerate}
    \item Compute a compact \emph{signature} $\Sigma(T_{new})$ as in the static case.
    \item Use blocking/indexing to retrieve a candidate set $\mathcal{C}_t \subseteq D_t$.
    \item Evaluate the multi-metric distance $\delta(T_{new}, T_i)$ for all $T_i \in \mathcal{C}_t$.
    \item Flag $T_{new}$ as a duplicate if $\delta(T_{new}, T_i) < \varepsilon$ for any $T_i$.
    \item Insert $T_{new}$ into $D_{t+1}$ with updated indices.
\end{enumerate}

\paragraph{Sliding Window or Bounded Memory Approach.}
For high velocity streams, storing all historical records is impractical. A
sliding window of the most recent $W$ records ensures that duplicate detection
remains computationally tractable:
\[
D_t = \{ T_{t-W+1}, T_{t-W+2}, \dots, T_t \}.
\]
Older records outside the window can be discarded or archived. Window based
duplicate detection is particularly effective when duplicates are temporally
localized.

\paragraph{Incremental Indexing and Pruning.}
Streaming duplicate detection benefits from indices that support insertion and
range queries:
\begin{itemize}
    \item \textbf{Signature insertion:} When $T_{new}$ arrives, its signature
    $\Sigma(T_{new})$ is immediately inserted into the index.
    \item \textbf{Candidate retrieval:} Only records within the same block
    or with similar signatures are queried.
    \item \textbf{Dynamic pruning:} Partial multi-metric computations can
    terminate early if a lower bound threshold is exceeded.
\end{itemize}

\paragraph{Complexity Analysis.}
Let $b$ denote the average block size. For each incoming record, the
streaming detection requires:
\[
O(b m) \quad \text{time per record},
\]
which is near constant when $b$ is bounded and $m$ is moderate. Memory
requirement is $O(W m)$ with a sliding window of size $W$, making the method
scalable for large and continuous datasets.

\paragraph{Analytical Example.}
Consider a stream of transactional records
\[
T_1, T_2, T_3, \ldots
\]
defined over the finite attribute universe
\[
X = \{a,b,c\}.
\]
Each record $T_i$ is represented as a multiset through its
multiplicity function
\[
C_{T_i} : X \to \mathbb{N},
\]
where $C_{T_i}(x)$ denotes the number of occurrences of
attribute $x$ in $T_i$.

Assume that duplicate detection is performed using a
multiplicity aware multi-metric
\[
\delta(T_i, T_j)
= \sum_{x \in A}
d\!\left(
P_a^{C_{T_i}(x)},\,
P_a^{C_{T_j}(x)}
\right).
\]

and let the duplication threshold be fixed as $\varepsilon = R^2_0$.

\medskip
\noindent
\textbf{Step 1: Arrival of $T_1$.}
The first record is
\[
T_1 = \{a^1,b^2,c^0\}.
\]
Since the index is initially empty, $T_1$ is inserted directly
into the index and serves as a reference record for future
comparisons.

\medskip
\noindent
\textbf{Step 2: Arrival of $T_2$.}
The second record is
\[
T_2 = \{a^1,b^1,c^0\}.
\]
Because $T_2$ shares the same support pattern as $T_1$,
it is placed in the same block and compared with $T_1$.
The multi-metric distance is computed as

\[
\begin{aligned}
\delta(T_1,T_2)
&= R^{|1-1|}_{0} + R^{|2-1|}_{0} + R^{|0-0|}_{0} \\
&= R^{0}_{0} + R^{1}_{0} + R^{0}_{0} \\
&= R^{1}_{0}.
\end{aligned}
\]

Since $\delta(T_1,T_2) < \varepsilon$, the record $T_2$ is
identified as a duplicate of $T_1$ and is accordingly flagged
and inserted into the index.

\medskip
\noindent
\textbf{Step 3: Arrival of $T_3$.}
The next record in the stream is
\[
T_3 = \{a^0,b^2,c^1\}.
\]
The support of $T_3$ differs from that of $T_1$ and $T_2$,
so blocking places $T_3$ in a separate group. For completeness,
consider the distance to $T_1$:

\[
\begin{aligned}
\delta(T_3,T_1)
&= R^{|0-1|}_{0} + R^{|2-2|}_{0} + R^{|1-0|}_{0} \\
&= R^{1}_{0} + R^{0}_{0} + R^{1}_{0} \\
&= R^{3}_{0}.
\end{aligned}
\]

Because $\delta(T_3,T_1) > \varepsilon$, $T_3$ is not classified
as a duplicate and is retained as a distinct record.

\medskip
\noindent
\textbf{Interpretation.}
This example illustrates how multi-metric based duplicate
detection operates efficiently in an online setting. By combining
incremental distance computation with blocking and threshold based
pruning, the method accurately detects duplicates while avoiding
unnecessary comparisons, even in a streaming environment.

\section{Conclusion}

Multi-metric spaces offer a mathematically robust extension of classical metric spaces
for frequency aware data analysis. Their application to duplicate detection demonstrates
clear advantages in preserving semantic distinctions arising from multiplicity.
Future work includes scalable algorithms and probabilistic extensions.

\end{document}